\documentclass[12pt]{article}
\usepackage{bbm}
\usepackage{fullpage}
\usepackage{graphicx}
\newtheorem{theorem}{Theorem} 

\newtheorem{lemma}[theorem]{Lemma}
\newtheorem{claim}[theorem]{Claim}

\newtheorem{definition}[theorem]{Definition}

\def\eod{\vrule height 10pt width 9pt depth 0pt}

\newenvironment{proof}{\noindent {\bf Proof:} \hspace{.4em}}
                      {\hspace*{\fill}{\eod}}

\newcommand{\ve}{\varepsilon}
\newcommand{\de}{\delta}
\newcommand{\al}{\alpha}
\newcommand{\be}{\beta}
\newcommand{\ga}{\gamma}
\newcommand{\ti}{\tilde}
\newcommand{\beq}{\begin{equation}}
\newcommand{\eeq}{\end{equation}}

\newcommand{\eref}[1]{(\ref{#1})}
\newcommand{\CC}{\mathbbm{C}}
\newcommand{\RR}{\mathbbm{R}}
\newcommand{\ra}{\rightarrow}
\newcommand{\lar}{\leftarrow}

\begin{document}

\centerline{\bf \Large Parabolic Julia Sets are Polynomial 
Time Computable}

\bigskip 

\centerline{Mark Braverman$^1$}
\centerline{Department of Computer Science}
\centerline{University of Toronto}
\bigskip
\centerline{\today}
\bigskip


\bigskip
\footnotetext[1]{Research is partially supported by an NSERC postgraduate scholarship.}

\begin{abstract}

In this paper we prove that parabolic Julia sets of 
rational functions are locally computable 
in polynomial time. 
\end{abstract}

\section{Introduction}

In the present paper we consider the complexity of generating 
precise images of Julia sets with {\em parabolic} orbits. 
It has been independently proved in \cite{Brv} and \cite{Ret} that 
{\em hyperbolic} Julia sets can be computed in polynomial 
time. Neither of the two algorithms can be applied in the 
parabolic case. In fact, both algorithms often slow down 
significantly as the underlying polynomial approaches one 
with a parabolic point. A na\"ive generalization of these 
algorithms would yield exponential time algorithms in 
the parabolic case, which are useless when one is trying to 
produce meaningful pictures of the Julia set in question. 

The same problem has been highlighted in the comments on 
computer graphics by John Milnor in \cite{Milnor}, Appendix H. 
The example considered there is for the polynomial $p(z) =
z+z^4$. It has a parabolic fixed point at $z=0$. 
Consider a point $z=\ve \approx 1/1000$. Suppose we are
trying to determine whether $\ve$ is in the Julia set or
not by iterating it, and observing whether its orbit escapes to 
$\infty$, or converges to $0$. In fact, such a $z$ would 
always escape to $\infty$, but it is not hard to see that 
this process would take $1/3 \ve^3 \approx 300,000,000$ iterations 
for $z$ to escape the ball of radius $2$ around $0$. Thus, we would need to 
follow the orbit for $\approx 300,000,000$ iterations before 
concluding that it converges to $\infty$. If we zoom-in a 
little and set $z=\ve \approx 1/100,000$, we would need $\approx 
3 \cdot 10^{14}$ iterations to trace $z$, which is 
computationally impractical. 

Due to the effect highlighted above, most computer programs 
plotting Julia sets include all the points that diverge slowly from 
the parabolic orbit in the Julia set. 

The algorithm we present here is not uniform, i.e. it requires
a special program for each specific parabolic Julia set. The 
running time of the algorithm is $C_r n^c$, where the 
constant $C_r$ depends on the rational function $r$ but not on $n$, and 
$c$ is some constant. The algorithm can be made uniform in the
$r$, provided some basic combinatorial information about the 
parabolic points. I.e. one algorithm can compute all the parabolic 
sets, if it is provided with some basic information about the rational function. 
The constant $C_r$ in the running time can still vary strongly 
for different functions $r$. For example, it is reasonable to 
expect that $J_{r_1}$ for $r_1 (z) = z + z^2$ would take less
time to compute than $J_{r_2}$ for $r_2 (z) = e^{2 \pi i /17} z + z^{23}$.
We prove the following:

\begin{theorem}
There is an algorithm $A$ that given 
\begin{itemize}
\item
a rational function $r(z)$ such 
that every critical orbit of $r$ converges either to an attracting 
or a parabolic orbit; and 
\item 
some basic combinatorial information about the parabolic orbits of $r$;
\end{itemize}
produces an image of the Julia set $J_r$. $A$ takes time $C_r n^c$ to decide 
one pixel in $J_r$ with precision $2^{-n}$. Here $c$ is some small constant and 
$C_r$ depends on $r$ but not on $n$. 
\end{theorem}

After this work was completed, John Milnor has informed us that he 
has used an algorithm similar to ours to produce pictures of Julia sets
with parabolic points. In particular, some of the pictures in \cite{Milnor}
were created this way. 

The rest of the paper is organized as follows. In section \ref{compprelim} 
we give the necessary preliminaries on the complexity theory over the reals. 
In section \ref{secstrat} we outline the general strategy for computing 
parabolic Julia sets fast. Sections \ref{growth}, \ref{nthiter} and \ref{iterz}
provide the main tool for the algorithm -- computing a ``long" iteration 
near a parabolic point. Finally, in section \ref{thealg}, we present and
analyze the algorithm.

\medskip
\noindent
{\bf Acknowledgment.} The author wishes to thank Ilia 
Binder and Michael Yampolsky 
for their insights and encouragement during the preparation 
of this paper. 

\section{Complexity over $\RR$ -- preliminaries}
\label{compprelim}
In this section we provide some preliminaries on the notion 
of complexity for sets and functions over $\RR^n$, in particular 
$\RR^2$. More 
details can be found in \cite{Brat}, \cite{Brv05} and \cite{WeiBook}.

\subsection{Complexity of Sets in $\RR^2$}
\label{Complexity}

Intuitively, we say the computational complexity of 
a set $S$ is $t(n)$ if it takes time $t(n)$ to decide whether to draw a pixel of size $2^{-n}$ 
in the picture of $S$. To make this notion precise, we have 
to decide what are our expectations from a picture of $S$. First of all, 
we expect a good picture of $S$ to cover the whole set $S$. On the 
other hand, we expect every point of the picture to be close to some
point of $S$, otherwise the picture would have no descriptive power about 
$S$. Mathematically, we write these requirements as follows:

\begin{definition}
\label{def1}
A set $T$ is said to be a $2^{-n}$-picture of a bounded set $S$ if 

(i) $S \subset T$, and (ii) $T \subset B(S,2^{-n}) =
\{ x \in \RR^2~:~|x-s|<2^{-n}~\mbox{for some~} s\in S\}$.
\end{definition}

Definition \ref{def1} is also equivalent to approximating $S$ by $T$ in 
the {\em Hausdorff metric}, given by
$$
d_H (S,T) := \inf \{ r : S \subset B(T,r)~~\mbox{and}~~
T \subset B(S,r) \}.
$$

Suppose we are trying to generate a  picture of a set $S$ using 
a union of round pixels of radius $2^{-n}$ with centers at all 
the points of the form $\left( \frac{i}{2^{n}}, 
\frac{j}{2^{n}}\right)$, with $i$ and $j$ integers. In order to 
draw the picture, we have to decide for each pair $(i,j)$ whether to 
draw the pixel centered at $\left( \frac{i}{2^{n}},
\frac{j}{2^{n}}\right)$ or not. We want to draw the pixel if 
it intersects $S$ and to omit it if some neighborhood of the pixel 
does not intersect $S$. Formally, we want to compute a function 
\beq
\label{star1}
f_S ( n, i/2^{n}, j/2^{n}) = 
\left\{ 
\begin{array}{ll}
1, & B((i/2^{n}, j/2^{n}), 2^{-n})\cap S \neq \emptyset \\
0, & B((i/2^{n}, j/2^{n}), 2 \cdot 2^{-n})\cap S = \emptyset  \\
0~~\mbox{or}~~1,~~& \mbox{in all other cases}
\end{array}
\right. 
\eeq

\begin{figure}[!h]
\begin{center}
\includegraphics[angle=0, scale=0.7]{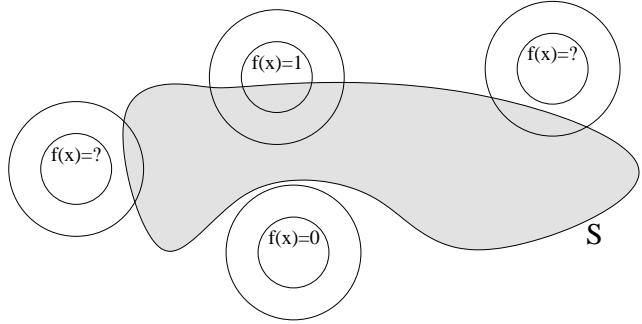}
\end{center}
\caption{Sample values of $f$. The radius of the inner circle is 
$2^{-n-2}$.}
\label{fig01}
\end{figure}

\begin{lemma}
The picture drawn according to $f_S(n,\bullet)$ is a $2^{-(n-2)}$-picture of 
$S$.
\end{lemma}

Here $\bullet$ stands for the different values of the parameters
$(i/2^{n}, j/2^{n})$.
The lemma illustrates the tight connection between the complexity of 
``drawing" the set $S$ and the complexity of computing $f$. We 
reflect this connection by defining the time complexity of  $S$ as follows.

\begin{definition}
\label{defcomp}
A bounded set $S$ is said to be computable in time $t(n)$ if there 
is 
a function $f(n,\bullet)$ satisfying \eref{star1} which runs in time 
$t(n)$. We say that $S$ is poly-time computable if there is a polynomial 
$p$, such that $S$ is computable in time $p(n)$. 
\end{definition} 

To see why this is the ``right" definition, suppose we are trying to 
draw a set $S$ on a computer screen which has a $1000\times 1000$ pixel
resolution. A $2^{-n}$-zoomed in picture of $S$ has $O(2^{2 n})$ pixels
of size $2^{-n}$, and thus would take time $O(t(n) \cdot 2^{-2 n})$ to 
compute. This quantity is exponential in $n$, even if $t(n)$ is bounded 
by a polynomial. But we are drawing $S$ on a finite-resolution display, 
and we will only need to draw $1000 \cdot 1000 = 10^6$ pixels. Hence 
the running time would be $O(10^6 \cdot t(n)) = O(t(n))$. This running 
time is polynomial in $n$ if and only if $t(n)$ is polynomial. 
Hence $t(n)$ reflects the `true' cost of zooming in. 

\subsection{Computing Julia Sets}
\label{compjuls}

There are uncountably many rational functions, but only countably many 
Turing Machines. Thus, we cannot expect to have  
a Turing Machine computing the Julia set $J_r$ for each rational $r(z)$. Instead, we assume that 
the coefficients of $r$ are given to the machine, and it is trying to 
produce a picture of $J_r$. The machine can access the coefficients
with an arbitrarily high (finite precision). It is charged $m$ time units
for querying a coefficient with precision $2^{-m}$. Hence if a machine 
computes $J_r$ with precision $2^{-n}$ in time polynomial $p(n)$, it 
will query the coefficients with precision at most $2^{-p(n)}$.

Another issue is whether the computation of a machine is uniform 
or non-uniform. A machine for computing $J_r$ is {\em non-uniform}, 
if it is designed specifically for this $r$. A machine is 
{\em uniform} on a set $S$ of rational functions, if it produces $J_r$ for 
all $r\in S$. One can view a non-uniform machine as a uniform machine
on the set $S = \{ r \}$. One of the properties of the computation 
model is that if $J_r$ is uniformly computable on $S$, then the 
function $J: r \mapsto J_r$ is continuous in the Hausdorff 
metric. In the case of a non-uniform computation, $S$ is a 
singleton, and thus we don't get any information from this 
statement. 

We first give a non-uniform algorithm for computing $J_r$. Then in section
\ref{unif} we argue that it can be made uniform for some large classes of 
parabolic Julia sets. The function $J: r \mapsto J_r$ is not continuous 
over all parabolic sets, and thus it cannot be uniformly 
computable on {\em all} parabolic functions $r$. See section 
\ref{unif} for more details.

\section{The Strategy}
\label{secstrat}

First we recall the strategy in the hyperbolic case, which is 
much easier to deal with. Suppose that  $r$ is a hyperbolic rational function.
Let $J_r$ denote its Julia set. Then $r$ is strictly expanding by some constant $c>1$ in the 
hyperbolic metric around $J_r$, and thus the escape rate of a point
$z \notin J_r$ near $J_r$ is exponential. In other words, if $d(z, J_r)>2^{-n}$, 
then after $O(n)$ steps the orbit of $z$ will be at $\Theta(1)$ distance
from $J_r$. This gives a natural poly-time algorithm for computing $J_r$:
iterate $z$ until it is possible to estimate the distance from $r^k(z)$ to $J_r$ 
using some coarse initial approximation to $J_r$. If such a $k = O(n)$ exists, 
use $d(r^k (z), J_r)$ and $|(r^k)'(z)|$ to estimate $d(z, J_r)$. If no such $k$ 
exists, we can be sure that initially $d(z, J_r)<2^{-n}$. 

\begin{figure}[!h]
\begin{center}
\includegraphics[angle=0, scale=0.47]{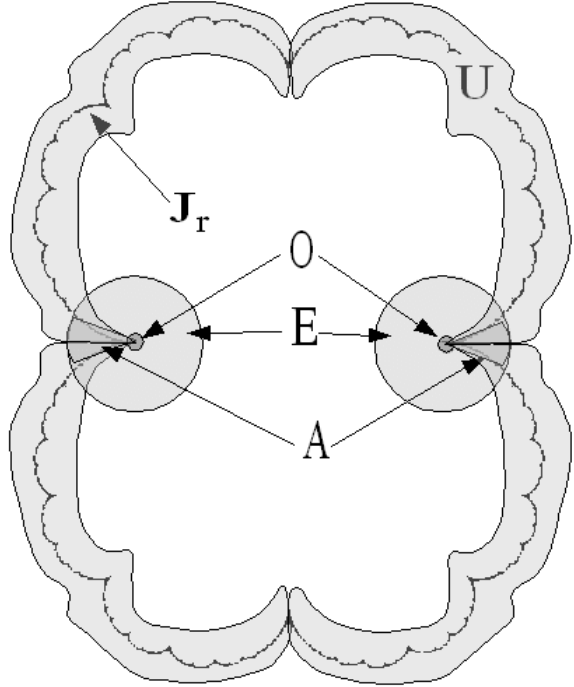}
\end{center}
\caption{A schematic image of the components in the algorithm.}
\label{stratf}
\end{figure}

We would like to employ a similar strategy here, in the parabolic case. The 
problem is that even though $r$ is still expanding in the hyperbolic metric
near $J_r$, the expansion is now  extremely slow near the parabolic point. 
For example, let $r(z) = z^2 + 1/4$, with the parabolic point $p = 1/2$. 
The picture of $J_r$ is presented of figure \ref{stratf}. If we set $z=1/2+2^{-n}$, 
it will take $O(2^n)$ steps before $z$ escapes the unit disk.

We solve this problem by approximating a ``long" iteration of $z$ in the neighborhood 
of a parabolic point fast. In the previous example ``long" would mean $O(2^{n})$.

On figure \ref{stratf}, we present the different regions which will appear in the 
algorithm. We list them below. 

\begin{itemize}
\item 
$J_r$ is the Julia set we are trying to compute. 
\item 
$U$ is some small fixed region around $J_r$. All the points in $U$ 
are much closer to $J_r$ than to the postcritical set. $U$ is 
bounded away from $J_r$, except for a finite number of parabolic and preparabolic
touching points. On figure \ref{stratf} the touching points are $1/2$ (the parabolic
point) and $-1/2$ (first-order preparabolic). 
\item
$E$ is the region around the parabolic points in which the ``long" iteration 
is applicable. We also include in $E$ preimages of this neighborhood around
the preparabolic touching points of $U$ and $J_r$. 
\item 
$A \subset U \cup E$ is a collection of small wedges around the repelling directions. 
These wedges contain the portions of $J_r$ in the neighborhood of the corresponding 
parabolic/preparabolic points. 
\item 
$O$ is a tiny neighborhood around the touching points of $U$ and $J_r$. If
the orbit of $z$ falls into $O$ we can be sure that $z$ is close to $J_r$ because 
all of $O$ is so close to $J_r$. 
\end{itemize}

Now the algorithm works exactly as the one in the hyperbolic case:

\begin{enumerate}
\item 
Iterate the orbit of $z$;
\item 
let $w=r^k(z)$ be the current iterate;
\item 
if $w \notin U \cup E$, we can estimate its distance from $J_r$ in $O(1)$ time; 
\item 
\label{sst4}
if $w \in U - E$, just make one step $w \lar r(w)$;
\item 
if $w \in O$, output ``$w$ close to $J_r$";
\item 
if $w \in E \cap A$ in the neighborhood of a preparabolic point,  just make one step $w \lar r(w)$;
\item 
if $w \in E- A$, we can estimate its distance from $J_r$ in $O(1)$ time;
\item 
\label{sst8}
if $w \in E \cap A$ near some parabolic point, apply linearly many ``long" iterations
to escape this region and get to step \ref{sst4};
\end{enumerate}

Step \ref{sst4} can only be executed linearly many times, since
$\partial U$ is bounded away from $J_r$ outside of $E$ and the 
expansion in the hyperbolic metric is bounded from below by some 
$c>1$ on $U-E$. Thus the entire computation takes at most a quadratic
number of steps to complete (at most linearly many executions of 
step \ref{sst8} between two executions of step \ref{sst4}). 

Of course, this is only a sketch, and we need more precise procedures taking
 into account the finite precision of the computation etc. (e.g. we cannot just
 check whether $w$ is in $A$ or not.) In the next sections we will develop
 the tools for performing the ``long" iteration near the parabolic points, 
 before formally presenting the algorithm.
\section{Controlling coefficient growth}
\label{growth}
The primary goal of this section is to prove the following 
lemma. 

\begin{lemma}
\label{growth:main}
Let $r \ge 1$ be some integer. Set
 $f(z)=z+z^{r+1}+z^{r+2}+z^{r+3} + \ldots$. Then there is 
 an explicit $\al$ such that the coefficients
 of the $n$-th iteration $f^n$ of $f$,
 $$
 f^n(z)= z +  a_r z^{r+1} + a_{r+1} z^{r+2} + a_{r+2} z^{r+3} + \ldots,
 $$
 satisfy
 $$
 a_k \le (\al n)^{k/r}.
 $$
 One can take $\al = 2 r^3$. 
\end{lemma}

We begin with a very simple proof in the case $r=1$. The general 
case is more involved. 

\medskip

\begin{proof} (in the case $r=1$.)
In this case $f(z) = z+z^2 + z^3+ \ldots = \frac{z}{1-z}$ (within 
the region of convergence).  It is easy to verify that
$$
f^n (z) = \frac{z}{1-n z}.
$$
Hence the coefficient of $z^{k+1}$ in $f^n$ is $n^{k}$, and lemma \ref{growth:main}
holds with $\al =1$. 
\end{proof}

\medskip

For the rest of the section we fix some $r \ge 2$, for which we are proving 
lemma \ref{growth:main}. We will prove the lemma by induction on $n$. 
It is obviously true for $n=1$ with any $\al \ge 1$.
Denote $\rho = (\al n)^{1/r}$. We will find a constant $\al$ that 
works later in the proof, $\al$ may depend on $r$ but not on $n$. 

If $g(z)$ and $h(z)$ are two power series with positive real coefficients, 
we say that $g$ is {\em dominated} by $h$, and write $g\ll h$ if all the coefficients
of $g$ are smaller or equal to the corresponding coefficient of $h$. 

We assume by the induction hypothesis that 
\beq
f^n(z) \ll z + \rho^{r} z^{r+1} + \rho^{r+1} z^{r+2} + \ldots.
\eeq

Denote $g(z) = z + r z^{r+1} + r z^{2 r +1 } + r z^{3 r+1} + \ldots$. We claim
the following. 

\begin{lemma}
\label{growth:l2}
Let $m \ge r$ be a given integer number. Then 
$$
f(z)^m -z^m \ll (1+z+z^2+\ldots+z^{r-1}) \cdot (g(z)^m-z^m).
$$
\end{lemma}

\begin{proof}
We show that the coefficient of $z^k$ on the left hand side 
is smaller or equal to the coefficient on the right hand side. 
Note that all the coefficients on the left for $k<m+r$ are $0$, and 
so we can assume that $k \ge m+r$. Write $k-m = r \cdot l + q$ 
with $0 \le q < r$, $l \ge 1$. We claim that the coefficient 
of $z^{k-q} = z^{r \cdot l + m}$ in $g(z)^m$ is bigger or
equal to the coefficient of $z^k$ in $f(z)^m$. 

To see this we create a one-to-one mapping from 
all the terms of degree $k$ in the expansion of 
$f(z)^m$ to the corresponding terms of degree $k-q$ in the expansion 
of $g(z)^k$ where we write
$$
g(z)=z^{1}_{(0)} + z^{r+1}_{(0)}+ z^{r+1}_{(1)} + \ldots + z^{r+1}_{(r-1)} + 
z^{2 r+1}_{(0)}+ z^{2 r+1}_{(1)} + \ldots + z^{2 r+1}_{(r-1)} + \ldots.
$$
Here $z_{(0)}, \ldots, z_{(r-1)}$ refer to different copies of 
the same $z$ (we separate the different copies to specify the 
one-to-one mapping). 

Suppose we are given a term $z^{a_1} z^{a_2} \ldots z^{a_m} = z^k$ in the 
expansion of $f(z)^m$. We write $a_i = 1 + r \cdot b_i + c_i$, $0 \le c_i <r$. 
Then we know that either $a_i-c_i = 1$ and $c_i=0$, or $a_i - c_i \ge r+1$. We associate the 
term 
\beq
\label{corr1}
z_{(c_1)}^{a_1-c_1} z_{(c_2)}^{a_2-c_2} \ldots z_{(c_{m-1})}^{a_{m-1}-c_{m-1}} 
z_{(0)}^{a_m + c_1 + c_2 + \ldots + c_{k-1} -q}.
\eeq
 By the construction $a_m + c_1 + c_2 + \ldots + c_{m-1} -q \equiv
 a_m + a_1 + a_2 + \ldots + a_{m-1} - (m-1)-q \equiv k - m - q + 1 \equiv 1 ~(mod~r)$.
$k \ge m + r$, so we will never need the term $z^m$ from $g(z)^m$. 
It is not hard to see that the correspondence is one-to-one, since the information in 
\eref{corr1} is sufficient to recover the values of $a_1, a_2, \ldots, a_m$. 

By considering the term $z^q \cdot g(z)^m$, we complete the proof. 
\end{proof}

\medskip

We are now ready to make the induction step in lemma \ref{growth:main}.
By the induction hypothesis and lemma \ref{growth:l2}, we have 
$$
f^{n+1}(z) \ll f(z) + \rho^r z^{r+1} + \rho^{r} (f(z)^{r+1}-z^{r+1})
 + \rho^{r+1} z^{r+2} + \rho^{r+1}  (f(z)^{r+2}-z^{r+2}) + \ldots \ll
$$
$$
f(z) + \rho^r z^{r+1}+\rho^{r+1} z^{r+2} + \rho^{r+2} z^{r+3} + \ldots + 
$$
$$
+(1 + z+ z^2 + \ldots + z^{r-1}) \cdot \left[ \rho^{r} (g(z)^{r+1}-z^{r+1}) + 
\rho^{r+1} (g(z)^{r+2}-z^{r+2}) + \ldots\right]
$$

Our goal is to bound the coefficient of $z^k$, $k \ge r+1$, in $f^{n+1}(z)= 
f^{n} (f(z))$. The contribution from $f(z)$ is always $1$. We consider 
the contribution from $g(z)^{m}-z^m$, $r+1 \le m \le k$. Write $k = m + r \cdot l + q$, 
$0 \le q < r$.  Then we must have $z^q$ in the product, and the coefficient 
is the coefficient of $z^{m + r\cdot l}$ in $\rho^{m-1} g(z)^m$, which is bounded by 
 $\rho^{m-1} {m+l-1 \choose l} r^l$. The contribution is nonzero only if $l >0$. 
 Thus, the coefficient is 
 bounded by 
 $$
1+\rho^{k-1}+ \sum_{q=0}^{r-1} \sum_{l=1}^{\left\lfloor \frac{k-q-r-1}{r}\right\rfloor}
 \rho^{k- rl -q -1} r^l {k-rl-q+l-1 \choose l} < 
 $$
 $$
 \rho^{k-1}+ r \sum_{l=1}^{\left\lfloor \frac{k-1}{r}\right\rfloor}
 \rho^{k- rl -1} r^{l} {k-rl+l-1 \choose l} < \rho^{k-1}+
  \sum_{l=1}^{\left\lfloor \frac{k-1}{r}\right\rfloor}
 \rho^{k- rl -1} r^{l+1} {k-1 \choose l}<
 $$
 $$
 \rho^{k-1}+ \rho^{k-1} \sum_{l=1}^{k-1}
 \rho^{- rl} r^{2 l} {k-1 \choose l} = 
 \rho^{k-1} \sum_{l=0}^{k-1}
 \rho^{- rl} r^{2 l} {k-1 \choose l} = \left(\rho+ \frac{r^2}{\rho^{r -1}}\right)^{k-1}. 
$$
To prove the lemma, we need the condition 
\beq
\label{growth:e3}
 \left(\rho+ \frac{r^2}{\rho^{r -1}}\right)^{k-1} \le 
 (\al (n+1))^{(k-1)/r}.
\eeq
Recall that $\rho = (\al n)^{1/r}$, hence we can rewrite \eref{growth:e3} as
$$
\left( (\al n)^{1/r} + \frac{r^2}{(\al n)^{(r-1)/r}} \right)^r  \le \al (n+1).
$$
We have 
$$
\left( (\al n)^{1/r} + \frac{r^2}{(\al n)^{(r-1)/r}} \right)^r = 
(\al n) \cdot \left( 1 + \frac{r^2}{\al n} \right)^r < (\al n) \cdot e^{\frac{r^3}{\al n}} <
(\al n) \cdot \left(1+ \frac{2 r^3}{\al n}\right),
$$
the last inequality holds whenever $\frac{r^3}{\al n}< \ln 2$. Finally, if we take $\al \ge 2 r^3$, then 
$$
(\al n) \cdot \left(1+ \frac{2 r^3}{\al n}\right) < (\al n) \cdot \left(1+ \frac{1}{n}
\right) = \al (n+1),
$$
as required.

\section{Computing the $n$-th iteration of $f$}
\label{nthiter}

Suppose that we are given a function $f$ presented 
as a power series, finite or infinite, $f(z)= z + a_2 z^2 + a_3 z^3 + \ldots$, 
$d \ge 2$. Denote for the $n$-th iteration $f^n$ of $f$
\beq
\label{nthiter:e1}
f^{n} (z) = z + a_2^{(n)} z^2 + a_3^{(n)} z^3 + \ldots
\eeq
The goal of this section is to show how to compute 
the values of $a_k^{(n)}$ with a given precision $2^{-l}$ fast  -- 
in time polynomial in $k$, $l$, and $\log n$. We will need this in order to interpolate 
``long" iterations of $f$ around the parabolic point ($0$ in this case). 

First, we show that if $f$ has a non-negative radius of 
convergence $R$, then we can assume that $|a_i| \le 1$ for all $i$, 
with a fairly small overhead. We know that $\sum_{i} a_i (R/2)^i$ converges. 
Hence there is a bound $B \ge 1$ such that $a_i (R/2)^i < B$ for all $i$. In 
other words, $|a_i| < B  \cdot (2/R)^i \le (2 B/R)^i/B$. In case $f$ is a 
rational function, it is easy to approximate the number $2 B/R$, or 
some power of two $C=2^c$, such that $|a_i| < C^{i-1}$ for all $i$. 
Conjugate $f$ by the map $z \mapsto C z$ to obtain $g(z) = C f(z/C)$.
Then $f(z) = g(C z)/C$, and $g^n (z)  = C f^n (z/C)$. The Taylor expansion 
of $g(z)$ is 
$$
g(z) = z + \frac{a_2}{C} z^2 + \frac{a_3}{C^2} z^3 + \ldots
$$
We see that all the coefficients of $g(z)$ do not exceed $1$ in
absolute value. The Taylor expansion of the $n$-th iteration 
of $g$ is 
$$
g^n (z) = z + \frac{a_2^{(n)}}{C} z^2 + \frac{a_3^{(n)}}{C^2} z^3 + \ldots
$$
 Thus, to compute $a_k^{(n)}$ with precision $2^{-l}$, we would 
 need to compute the coefficient of $z^k$ in $g^n$ with precision
 $2^{-(l+c\cdot k)}$. This is a linear overhead, and if we 
 can compute approximations for $g$ in time $poly(k,l,\log n)$, we
 will also be able to do it for $f$. From now on, we assume that 
 $|a_i| \le 1$ for all $i$. 
 
 \medskip

We prove the following lemma. 

\begin{lemma}
\label{nthiter:main}
Suppose $f(z) = z + a_2 z^2 + a_3 z^3 + \ldots$ is given 
by its power series. Then the coefficient $a_k^{(n)}$ as 
in \eref{nthiter:e1} can be presented by a polynomial 
of the form 
\beq
\label{nthiter:e2}
a_k^{(n)} = \al_0^{k} + \al_1^k n + \al_2^k n^2 + \ldots + \al_{k-1}^{k} n^{k-1},
\eeq
and the values of $\al_i^{j}$ for $j=2,3, \ldots, k$, $i=0, \ldots, k-1$, 
 can be computed with precision $2^{-s}$ in 
time polynomial in $k$ and $s$.
\end{lemma}

\begin{proof}
We prove the lemma by providing an iterative algorithm that 
computes the values of $\al_i^{k}$. In the same time we
prove that $a_k^{(n)}$ is indeed of the form as in equation
\eref{nthiter:e2}. On the $j$-th iteration we compute the
values of $\al_0^{j}$, $\al_1^{j}$, \ldots, $\al_{j-1}^{j}$.

Here is how to compute the $\al_{i}^j$ from the $\al_{i}^{j-1}$'s. We know 
that 
$$
f^n (z) = f^{n-1} (f(z)) = f(z) + a_2^{(n-1)} \cdot f(z)^2 +
a_3^{(n-1)} \cdot f(z)^3+ \ldots
$$
In order to compute $a_j^{(n)}$ we need to find the coefficient 
of $z^j$ in each of the terms $f(z)$, $f(z)^2$, $\ldots$, $f(z)^{j}$ 
(it is always $0$ in higher terms). All we have to do is to compute 
all the coefficients of $f(z), f(z)^2, \ldots, f(z)^j$ up to $z^j$
with precision $2^{-m}$ in time polynomial in $m$ and $j$.

We assume here that   the coefficients 
$a_2, a_3, \ldots, a_j$ are given as oracles. 
In case that $f$ is a rational function
$$
f(z) = \frac{z+ p(z) \cdot z^2 }{c - q(z) \cdot z},
$$ 
we know that $c \neq 0$, by the parabolicity of $f$ around $0$, hence 
$$
f(z) = \frac{1}{c} \cdot \frac{z+ p(z) \cdot z^2 }{1- q(z) \cdot z/c} = 
\frac{1}{c} \cdot \sum_{k=0}^{\infty} (z+p(z) z^2) \left( \frac{q(z) \cdot z}{c} \right)^k.
$$
The first $j$ coefficients of the expansion are now easily computed 
from this last formula.

The computation is done using a simple ``doubling" algorithm:
first compute $f(z)^2$, $f(z)^{2^2}$, $f(z)^{2^3}$, $\ldots$, $f(z)^{2^{\left\lfloor \log j \right\rfloor}}$,
and then compute the desired power of $f(z)$ as a combination of these. At each
multiplication we ``chop" all the terms of degree $j+1$ and above, hence
the entire computation is polynomial. All the coefficients at 
all times are bounded by $j^j$ in absolute value, hence the error 
is multiplied by at most $2 j^j = 2^{O(j \log j)}$ at each step, and 
we will need to do the operations with a precision of $2^{-(m + O(j \log^2 j))}$
in order for the final coefficients to be with precision $2^{-m}$. 

Note that in the case when $f$ is a finite degree polynomial, we can 
evaluate the coefficients of its powers using the multinomial formula,
and with no need for the numerical iterative computation described above.

Once we have the coefficient $c_j^{(i)}$ of $z^j$ in $f(z)^i$, we are
able to write
\beq
\label{nthiter:e3}
a_j^{(n)} = a_j + a_2^{(n-1)} c_j^{(2)} + a_3^{(n-1)} c_j^{(3)} + \ldots + 
a_j^{(n-1)} c_j^{(j)}.
\eeq
We already have all the parameters in \eref{nthiter:e3}, as numbers or explicit 
polynomials in $n$ of degree at most $j-2$, except for $a_j^{(n-1)}$. It is 
easy to see that $c_j^{(j)} = 1$. Thus, we obtain an explicit recurrence, that 
connects $a_j^{(n)}$ with $a_j^{(n-1)}$, and yields
$$
a_j^{(n)} = a_j + \sum_{i=1}^{n-1} \left( 
a_j + a_2^{(i)} c_j^{(2)} + a_3^{(i)} c_j^{(3)} + \ldots + 
a_{j-1}^{(i)} c_j^{(j-1)}\right).
$$
Thus, $a_j^{(n)}$ is given by a polynomial of degree at most $j-1$ in $n$. 
The coefficients can be computed very efficiently (see \cite{Concrete} for
more information on how to compute the sum $\sum_{i=1}^n i^j$). The precision 
bit loss in this process is also limited to $O(k \log^2 k)$ bits. The 
coefficients of the polynomial $a_k^{(n)}$ are precisely the information
we are looking for to complete the proof of the lemma. 
\end{proof}

\section{Computing a ``long" iteration}
\label{iterz}

We can now apply the results of sections \ref{growth} and  \ref{nthiter}
to prove the following lemma.

\begin{lemma}
\label{iterz:main}
Suppose $f(z) = z + a_r z^{r+1} + a_{r+1} z^{r+2} + \ldots$ is given 
by its power series with some 
positive radius of convergence $R$. Then there is an easily computable number
$C$ such that if $|z|<\frac{1}{m}<R$, we can compute the 
$\ell= \left\lfloor \frac{m^r}{C} \right\rfloor$-th iterate of $z$ and 
its derivative $\frac{d f^{\ell}}{d z}(z)$ with precision $2^{-s}$ in
time polynomial in $s$ and $\log m$. 

The loss of precision from $z$ to $f^{\ell}(z)$ 
can be bounded to a constant number of bits. The loss
of precision from $z$ to  $\frac{d f^{\ell}}{d z}(z)$ is 
$O(-\log |z|)$ bits. 
\end{lemma}

\begin {proof}
We begin similarly to the discussion in the beginning of 
section \ref{nthiter}. As before, it is easy to compute 
a power of two, $A = 2^{a}$ such that $|a_i|<A^{i-1}$
for all $i$. Again, let $g(z) = A f(z/A)$. Then all the 
coefficients of $g$ are bounded by $1$ in absolute value. 
Write 
$$
g^\ell (z) = z + b_r^{(\ell)} z^{r+1} + b_{r+1}^{(n)} z^{r+2} + \ldots
$$ 
$g$ is dominated by the series $z+z^{r+1}+z^{r+2}+\ldots$. 
Thus, using lemma \ref{growth:main} we conclude that $|b_k^{(\ell)}|<
(\al \ell)^{k/r}$ for some simple, computable $\al$. If we write 
\beq
\label{iterz:e1}
f^\ell (z) = z + a_r^{(\ell)} z^{r+1} + a_{r+1}^{(\ell)} z^{r+2} + \ldots = 
\frac{1}{A} g^\ell (A z), 
\eeq
we see that $a_k^{(\ell)} = A^k b_k^{(\ell)}$, and $|a_k^{(\ell)}|<(\al A^r \ell)^{k/r}$.
Considering  that $\ell = \left\lfloor \frac{m^r}{C} \right\rfloor$ and $|z|<\frac{1}{m}$, we
obtain 
$$
|a_k^{(\ell)} z^{k+1}|< \left(\frac{\al A^r m^r}{C m^{r+1}}\right)^{k/r}< \left(\frac{\al^{1/r} A}{C^{1/r}}\right)^k.
$$
Choose $C> 2^r \al A^r$. Then $|a_k^{(\ell)} z^{k+1}|<2^{-k}$, and it suffices to consider the
first $s+2$ terms of the series \eref{iterz:e1} to obtain the desired iteration with 
a $2^{-s}$ precision (all later terms become negligible). 

All we have to do now is to compute $a_r^{(\ell)}$, $a_{r+1}^{(\ell)}$, $\ldots$, $a_{s+2}^{(\ell)}$ with 
precision $2^{-(s+\Theta(1))}$. We do it by computing their coefficients from \eref{nthiter:e2} with precision
$2^{-(s+O(s \log \ell))}$, which can be done in time polynomial in $s$ and $\log \ell = \Theta (\log m)$ by 
lemma \ref{nthiter:main}.

To compute the derivative of the $\ell$-th iteration, write
$$
\frac{d f^\ell}{d z} (z) = 1 + (r+1) a_r^{(\ell)} z^{r} + (r+2) a_{r+1}^{(\ell)} z^{r+1} + \ldots 
$$
then 
$$
|(k+1) a_k^{(\ell)} z^{k}| < (k+1) \left(\frac{\al A^r m^r}{C m^{r}}\right)^{k/r} < (k+1) \cdot 2^{-k} < 2^{-k/2},
$$
for sufficiently large $k$. Hence  it suffices to consider the
first $2 s+2$ terms of the series \eref{iterz:e1} to obtain the desired iteration with 
a $2^{-s}$ precision (all later terms become negligible). 
We compute $a_r^{(\ell)}$, $a_{r+1}^{(\ell)}$, $\ldots$, $a_{2 s+2}^{(\ell)}$ with 
precision $2^{-(s+\Theta(1))}$, which again can be done in time polynomial in $s$ and $\log \ell = \Theta (\log m)$ by 
lemma \ref{nthiter:main}.

The loss of precision can be kept to a constant number 
of bits for $z$ by the constant bound we have on 
the first derivative of $f^{\ell}(z)$ around $z$. The 
second derivative is bounded by $O\left({1}/{|z|}\right)$ around 
$z$, and the precision loss can be kept to $O(-\log |z|)$ bits, 
which is fine, as long as $|z|$ is not too small. 
\end{proof}

\section{Computing parabolic Julia sets in polynomial time}
\label{thealg}

In this section we put the pieces together to give
a poly-time algorithm for computing parabolic 
Julia sets. For the rest of the section fix $r(x)$ 
to be a rational function on $\CC$, and denote its
Julia set by $J_r$. We consider $J_r$ first as a 
subset of the Riemann sphere $\hat{\CC}$. Using 
a stereographic projection $\pi:\hat{\CC}-\{\infty\} \ra \CC$, 
 for any compact $C \subset \hat{\CC}$ such that $\infty \notin C$, it is easy to see that 
computing $J_r$ on $C$ is exactly as hard as computing 
$\pi (J_r)$ in $\CC$. Hence, if $\infty$ is not in $J_r$, 
it suffices to compute it in some bounded region $B(0, R)$ in $\CC$. 

If $\infty \in J_r$, then it is obviously impossible to 
``draw" it on the plane. We can still ``draw" it on the 
Riemann sphere, and hence on any bounded region of $\CC$. 
We take a M\"{o}bius transformation $T$ such that for the
conjugation $r' = T \circ r \circ T^{-1}$, $J_{r'}$ is 
obtained from $J_r$ by a rotation of the Riemann sphere. 
In this way, drawing $J_{r'}$ on $\hat{\CC}$ is as easy
(or as difficult) as drawing $J_r$. We can choose 
$T$ so that $\infty \notin J_{r'}$, and then it suffices 
to draw $J_{r'}$ on some bounded region of $\CC$. 

From now on, we assume that $\infty \notin J_r$, and 
that we have some constant $B$ such that $J_r \subset B(0,B)\subset \CC$. 
We are trying to compute $J_r$ on this bounded region.

\subsection{Preliminaries -- the nonuniform information we will need}
\label{nonunif}

Below we list the information the algorithm will use to 
compute $J_r$ efficiently with an arbitrarily high precision. 
We will need  the following ingredients:
\begin{enumerate}
\item 
A list of periods $v_1, v_2, \ldots, v_k$ for all the parabolic 
orbits. Consider the iteration $r^v (z)$ of $r(z)$ for $v = LCM (v_1, v_2, \ldots, v_k)$.
It  has only simple parabolic points (no orbits). 
From now on, we replace $r(z)$ with $r^v (z)$. We can do it, because 
$J_{r^v} = J_r$ for all $v$. We can multiply  $v$ by some other factor, 
so that the derivative $\frac{d r^v(z)}{dz}$ at each parabolic point is
$1$, and not some other root of unity. 
\item 
Information that would allow us to identify the parabolic 
points, and information about them. For each parabolic point $p$  
we would like to know an approximation $q$ for $p$ that would 
allow us to compute $p$ in poly-time using Newton's method. Near $p$, 
$r (z)$ can be written as 
$$
r  (z) =  p+  (z-p) + \al_{u} (z-p)^{u+1} + \al_{u+1} (z-p)^{u+2} + \ldots
$$
for some integer  $u$. We would like to know this number. 
\item 
\label{pr3}
An open set $U$ such that 
\begin{itemize}
	\item 
	$U \supset V := r^{-1} (U)$, 
	\item 
	all the parabolic points are in $\partial U$, 
	\item
	$J_r \subset \overline{U}$,
	\item
	all the critical (and hence also all the postcritical) points of $r(z)$ lie 
	outside $U$, 
	\item 
	all the poles and their neighborhoods lie outside $U$, 
	\item 
	moreover, if we denote the postcritical set by $P_r$, then for any $u \in r(U)$, 
	$d(u, P_r) \ge 32 \cdot d(u, J_r)$, and  
	\item 
	outside any $\ve$-neighborhood of the parabolic points and a finite 
	number of their pre-images, the distance between 
	$\partial U$ and $\partial V$ is bounded from below by some positive $\de$.
\end{itemize}
$U$ is given in the form $U = r^{-u} (\ti{U})$, where $\ti{U}$ is some explicit 
semi-algebraic set. Thus queries about membership in  $U$ and $V$ can be computed efficiently with an arbitrarily
high precision, at least outside some small region around the parabolic points and their 
preimages up to order $u$. We will show that such a $U$ exists, and 
how to compute it from some basic combinatorial information in section \ref{unif}.
\item 
\label{pr4}
For each parabolic point $p$, there is a small neighborhood $E_p$ of $p$ 
in which lemma \ref{iterz:main} applies for computing a long iteration 
of $r$. We would like to have two  sets $E_1$ and $E_2$ around the parabolic points
and their pre-images of order up to $u$ such that
\begin{itemize}
	\item 
	$E_1 \subset E_2$,
	\item 
	for a given point $z$, it takes constant time to decide whether $z\in E_1$, or 
	$z \notin E_2$, 
	\item 
	for each $z \in E_2$, there is a parabolic point $p$ such that $w = r^u (z) \in E_p$, 
	\item
	$\partial U \cap \partial V \subset E_1$, and
	\item
	we have a positive $d$ such that for any two points $x_1 \in V$, $x_2 \notin U$ outside
	of $E_1$, $|x_1-x_2|>d$.
\end{itemize}

\item 
\label{pr5n}
The set $\partial U \cap \partial V$ consists of pre-parabolic points $q$, i.e.
points such that $r^{u} (q)$ is parabolic for some fixed $q$. The repelling 
directions and their pre-images belong to $V$. There is an angle 
$\al$ such that all the points in $E_2$ that form  an angle $<\al$ with
		one of the repelling directions, or their preimages belong to $V$.
		If necessary, we can make $E_2$ smaller. We denote the subset in 
		$E_2$ of points that make an angle of $<\al/2$ with a repelling direction 
		or its pre-image by $A_1$, and the  points that make an angle of $<\al$ with a repelling direction 
		or its pre-image by $A_2$. We can choose $\al$ as small as we want. 
		We have the following properties.
\begin{itemize}
		\item 
		$A_1 \subset A_2 \subset  V$,
		\item
		if $z$ is given within an error of $<|z-q|^2$, near a 
		pre-parabolic point $q \in r^{-u} (p)$, we can tell if $z\in A_1$ or 
		$z \notin A_2$,
		\item 
		for any $p$, and  for any  $w \in A_2 \cap E_p$, 	$|r'(w)|>1$. This 
		is true for a sufficiently small $\al$. 
\end{itemize}

\item 
\label{pr5}
Consider the Poincar\'e metric defined on the hyperbolic set $U$.
Denote its density by $d_U$.  We have the following theorem,
known as Pick's theorem (see \cite{Milnor} for a proof).
\begin{theorem}
\label{Pick}
{\bf (Theorem of Pick)} Let $S$ and $T$ be two hyperbolic subsets 
of $\CC$. 
If $f:S \rightarrow T$ is a holomorphic map, then exactly one of the 
following three statements is valid:

\begin{enumerate}
\item
 $f$ is a conformal isomorphism from $S$ onto $T$, and maps $S$ with 
its Poincar\'{e} metric isometrically to $T$ with its Poincar\'{e} metric. 
\item 
 $f$ is a covering map but is not one-to-one. In this case, it is locally 
but not globally a Poincar\'{e} isometry. Every smooth path $P:[0,1] \rightarrow
S$ of arclength $l$ in $S$ maps to a smooth path $f \circ P$ of the same 
length $l$ in $T$. 
\item 
 In all other cases, $f$ is a strict contraction with respect to the
Poincar\'{e} metrics on the image and preimage.
\end{enumerate}
\end{theorem}

Let $d_V$ be the density of the Poincar\'e metric defined on $V$. By the construction, 
$V$ contains no critical points, and so $r: V \ra U$ is a covering map
and by theorem \ref{Pick} it is a local isometry. That is, for any 
$z \in V$, 
\beq
\label{thealg:e1}
d_V(z)  = d_U (r(z)) \cdot |r'(z)|.
\eeq
 On the other hand, 
the embedding $\iota : V \hookrightarrow U$ is not a covering map, 
hence it is strictly contracting in the Poincar\'e metric. Thus for any
 $z \in V$ we have $d_V(z)>d_U(z)$. Together with \eref{thealg:e1}, this 
 implies
 \beq
 \label{thealg:e2}
 d_U (r(z)) \cdot |r'(z)| = d_V (z)> d_U (z)
 \eeq
 for all $z \in V$. In particular, if we consider only $z \in V-E_1$, 
 then $z$ is in some compact domain bounded away from the boundary of
 $U$, hence the ratio $d_V(z)/d_U(z)$ is always positive (maybe $\infty$), 
 and it has a minimum $c>1$. We would like to have this $c$ as part
 of the nonuniform information. With this $c$ we have $d_V (z) > c \cdot d_U (z)$ 
 for all $z \in V-E_1$, and \eref{thealg:e2} becomes
 \beq
 \label{thealg:e3}
 d_U (r(z)) \cdot |r'(z)| = d_V (z)> c \cdot d_U (z)
 \eeq
 for all $z \in V-E_1$. Moreover, we can choose a slightly 
 smaller $c$ such that \eref{thealg:e3} holds for any
 point $z$ on any path $p$ from $w \in V-E_1$ to $J_r$ such that 
 $L_{d_U}(p)\le 2 d_{d_U} (w, J_r)$. This is true since the lengths 
 of such paths can be uniformly bounded, and thus it
  cannot get too close to the points of $\partial U \cap 
 \partial V$. 
 We would like to have the value of $c$ (or some rational 
 estimate $c_0$, $1< c_0 \le c$). 
 
 \item 
 \label{pr6}
 Since the postcritical points are outside $U$, and the parabolic 
 points cannot be critical, we can have a constant $d>0$ such 
 that $|r'(z)|>d$ for all $z \in r(U) \cup E_2$ (if necessary, 
 we can choose a smaller $E_2$. 
 
 \item 
 \label{pr8}
 Finally, we need an efficient procedure to estimate the distance from $J_r$ 
 for all points that are not too close to it. More specifically, for any point $z$
 outside of $V \cup E_1$, there is an ``estimator" that provides the distance
 $d(z, J_r)$ within a multiplicative error factor of $2$. This can be done since
 $J_r \subset V$ and the distance $d( \partial (V \cup E_1), J_r)$ is bounded 
 from below by a constant. Hence a fixed-precision image of $J_r$ suffices
 to make such an estimation. 
 
 The situation is somewhat different if $z \in E_2$. In this case, we 
 know that either $z$ or $r^k (z)$ for some bounded $k$ is close 
 to a parabolic point $p$. We know that in some small neighborhood 
 of $p$, $J_r$ looks like $\ell$ lines at angle $2 \pi / \ell$ from 
 each other leaving $p$ (see \cite{Milnor}, Chapter 10 for more details). Denote this set 
 by $L_p$. In general, if $z$ is very close to the Julia set, 
 $d(z, L_p)$ can be (multiplicatively) very different from 
 $d(z, J_r)$. However, if we stay away from the repelling directions
 (and hence from $L_p$ and $J_r$), these two quantities are actually
 similar. More precisely, for any $\al$ (and in particular for 
 $\al$ mentioned in the definition of $A_1$),
 in a small neighborhood of $p$, $d(z, L_p)$ is within a factor 
 of $2$ from $d(z, J_r)$ for every $z$ that
 has an angle of at least $\al/2$ with each of 
 the repelling directions. The same holds for the first $u$ 
 pre-images of the parabolic points. 
 We can take $E_2$ (and hence $A_1$) to be sufficiently small
 such that this property holds within $E_2-A_1$. 
\end{enumerate}

\subsection{Algorithm outline and analysis}

The goal of the algorithm is to compute a function 
from the family
$$
f(z,n) = \left\{ 
\begin{array}{ll}
1, & \mbox{if } d(z,J_r) < 2^{-n} \\
0, & \mbox{if } d(z,J_r) > 256 \cdot 2^{-n} \\
0 \mbox{ or } 1, & \mbox{otherwise}
\end{array}
\right.
$$
 To do this, we estimate 
$d (z, J_r)$ up to a multiplicative constant, assuming that
$d (z, J_r)>2^{-n}$. If the assumption does not hold, the 
algorithm always outputs $1$ (or successfully estimates the 
distance). 

The algorithm outline is as follows. 

\begin{enumerate}
\item 
$w \lar z$; 
$steps \lar 0$;
\item 
if $w \notin V \cup E_1$, output $0$;
\item
\label{st3}
estimate the maximum number $N = O(n)$ of steps 
outside $E_1$ we would need;
\item
iterate the point $w$ as follows:
\item 
set derivative counter $D \lar 1$; the cumulative derivative estimation should be bounded 
between the derivative and twice the derivative at all steps;
\item 
\label{st6}
if $w \notin V \cup E_1$:
\begin{enumerate}
	\item 
	estimate $e$, a $2$-approximation of $d(w, J_r)$: $e \le d(w, J_r) \le 2 e$; 
	\item 
	output $1$ if $\frac{e}{D} \le 8 \cdot 2^{-n}$, and $0$ if $\frac{e}{D} \ge 16 \cdot 2^{-n}$;
\end{enumerate}
\item 
\label{st7}
if $w \in U-E_1$:
\begin{enumerate}
	\item
	$D \lar D \cdot |r'(w)|$
\item
 $w \lar r(w)$;
 \item 
 $steps \lar steps + 1$;
 \item
 if $steps > N$, output $1$;
\end{enumerate}
\item 
\label{st8}
if $w \in E_2$ in the region of some preparabolic  $q\in r^{-u}(p)$, and $|w-q|<2^{-\be n}$ ($\be$ -- a constant
to be determined):
\begin{enumerate}
	\item
	output $1$;
\end{enumerate}
\item 
\label{st9}
if $w \in E_2 \cap A_2$, and it is not in the  neighborhood of any parabolic point:
\begin{enumerate}
	\item
	$D \lar D \cdot |r'(w)|$
\item
 $w \lar r(w)$;
\end{enumerate}
\item
\label{st10}
if $w \in E_2 - A_1$:
\begin{enumerate}
	\item 
	estimate $e$, a $2$-approximation of $d(w, J_r)$: $e \le d(w, J_r) \le 2 e$; 
	\item 
	output $1$ if $\frac{e}{D} \le 8 \cdot 2^{-n}$, and $0$ if $\frac{e}{D} \ge 16 \cdot 2^{-n}$;
\end{enumerate}
\item
\label{st11}
if $w \in E_p \cap A_2$ for some $p$:
\begin{enumerate}
\item
make a long iteration $y = r^v (w)$;
\item 
if $y \in E_p$, but escapes $A_2$:
\begin{enumerate}
\item
 do binary search to find the smallest 
$u<v$ such that $r^u (w)$ is in $A_2 - A_1$;
\item 
$w\lar r^u (w)$
\item
 go to step \ref{st10};
\end{enumerate}
\item 
else, 
$w \lar r^v (w)$;
\item
$D \lar D \cdot \left| \frac{d r^v (w)}{d w} \right|$;
\end{enumerate}
\end{enumerate}

The algorithm performs the operations with $O(n^2)$ bits
of precision. First note that steps 
\ref{st6}--\ref{st11} cover all the possibilities for $w$. 
If two or more of the possibilities intersect, it does not 
matter which one to choose. 

We first show that 
\begin{claim}
\label{cl3}
Step \ref{st3} in the algorithm is possible. 
\end{claim}

\begin{proof}
Let $z$ be some point in $V$ outside $E_1$. Let $p$ be 
the shortest path in the Poincar\'e metric $d_U$ from $p(0) =f(z)$ to $p(1) \in J_r$. $r: V \ra U$ 
is a covering map, and $p$ can be raised to a path $\ti{p}$ 
in $V$ such that $\ti{p} (0) = z$ and $\ti{p} (1) \in J_r$ (by
the invariance of $J_r$ under $r$). There are two possibilities:
\begin{enumerate}
\item 
$L_{d_U} (\ti {p}) \ge 2 d_{d_U} (z, J_r)$. We know that 
$r$ expands the Poincar\'e metric $d_U$, and hence
$$
d_{d_U} (r(z), J_r) = L_{d_U} ({p}) \ge L_{d_U} (\ti {p}) \ge 2 d_{d_U} (z, J_r) \ge
c \cdot d_{d_U} (z, J_r).
$$
\item
Otherwise, by property \ref{pr5} of $c$ from section \ref{nonunif}, $r$ is 
expanding by a factor of $c$ along the entire path $\ti{p}$. Hence
\beq
\label{expan1}
d_{d_U} (r(z), J_r) = L_{d_U} ({p}) \ge c \cdot L_{d_U} (\ti {p})  \ge
c \cdot d_{d_U} (z, J_r).
\eeq
\end{enumerate}
This shows that every step \ref{st7} multiplies the Poincar\'e 
distance between $w$ and $J_r$ by a factor of at least $c$. 
Other steps do not decrease it. If initially
the Euclidean distance $d(w, J_r)$ is at least $2^{-n}$, then the 
Poincar\'e distance is at least $C \cdot 2^{-n}$ for some constant
$C$. For every point in $V-E_1$ this distance is bounded from above, 
hence it will take at most $\log_c C^{-1} \cdot 2^{n} = O(n)$ steps 
\ref{st7} for the orbit of $z$ to escape.
\end{proof}

\begin{claim}
\label{clinU}
At any stage of the algorithm, $w \in r(U)$. 
\end{claim}

\begin{proof}
If on the first iteration we do not exit on step \ref{st8} or \ref{st10}, 
then we must have $w\in A_2 \cup V = V \subset U$ before the first iteration. 
From here, we prove the claim by induction. 

Suppose the algorithm is running after $i$ iterations. Denote the current
value of $w$ by $w'$ and the value after the iteration by $w''$ ($w''=w'$ if
the algorithm terminates). We assume that $w' \in r(U)$.
If the algorithm executes steps \ref{st7} or \ref{st9}, then by the conditions 
$w' \in U$, and $w'' = r(w') \in r(U)$. Quitting on steps \ref{st8} and 
\ref{st10} does not affect the value of $w$. Step \ref{st11} only runs
to keep $w''$ in $A_2 \subset V$. Thus $w''\in r(U)$ in this case as well. 
\end{proof}

\medskip 

The following is a classical theorem in complex analysis. 

\begin{theorem}{\bf Koebe's $1/4$ Theorem}
Suppose $\phi: (S_1, s_1) \ra (S_2, s_2)$ is a holomorphic bijection between 
two simply connected subsets of $\CC$, $S_1$ and $S_2$. Let $r_1$ be the 
inner radius of $S_1$ around $s_1$, and $r_2$ be the inner radius of 
$S_2$ around $s_2$. Then the following inequality holds:
$$
r_2 \ge \frac{1}{4} r_1 \cdot |\phi'(s_1)|.
$$
\end{theorem}

We apply Koebe's theorem to prove the following lemma.

\begin{lemma}
\label{Kcol}
Suppose that $z \in r(U)$ and the algorithm terminates with $w = r^{k} (z)$,
$D \le \left|\frac{d r^k (z)} {d z} \right| \le 2 D$, and 
$e \le d(w,J_r)\le 2 e$ for some distance estimate $e$. 
Then 
\beq
\label{k1}
\frac{e}{8 D} \le d (z, J_r) \le \frac {16 e}{D}.
\eeq
\end{lemma}

\begin{proof}
Denote $a = \left|\frac{d r^k (z)} {d z} \right|$ and  
$s = d(w, J_r)$. Then since $w \in r(U)$ the distance 
$d(w, P_r)$ from the postcritical set is at least $32 s$. 
We consider orbit $z$, $r(z)$, $r^2 (z)$, $\ldots$, $r^k (z)=w$.
Let $S_k = B(w, 32 s)$. Consider the preimage $S_{k-1}$ of
$S_k$ under the branch of $r$ that takes $r^{k-1}(z)$ to $r^k(z)$. 
It is uniquely defined since $S_k$ contains no postcritical points. 
The mapping $r: S_{k-1} \ra S_k$ is a one-to-one conformal 
mapping. We can continue this process to obtain a one-to-one 
conformal branch $r^k: (S_0, z) \ra (S_k, w)$. By Koebe's theorem 
the image of $B(w, s)$ under the inverse of this mapping 
must contain a ball of radius at least $\frac{s}{4 a}$ around $z$.
By the invariance of $J_r$, this ball contains no points from $J_r$.
Hence  
$$
d(z, J_r) \ge \frac{s}{4 a} \ge \frac{e}{4 \cdot 2 D} = \frac{e}{8 D}. 
$$
Also by Koebe's theorem, $S_0$ contains the ball $B_1$ of radius 
 $\frac{32 s}{4 a} = \frac{8 s}{a}$ around $z$. The image $r^{k}(B_1)$
must contain a ball of radius $ \frac{8 s}{a} \cdot \frac{a}{4} = 2 s$ around
$w$. Hence $r^{k} (B_1)$, and also $B_1$ contain points from $J_r$. 
So 
$$
d(z, J_r) \le r(B_1) = \frac{8 s}{a} \le \frac{8 \cdot 2 e}{D} = \frac{16 e}{D}.
$$
\end{proof}

\begin{claim}
If the algorithm terminates at step \ref{st6} or \ref{st10}, 
it outputs a valid answer. 
\end{claim}

\begin{proof}
The variables $e$ and $D$ in these cases satisfy the conditions
of lemma \ref{Kcol}. If $d(z, J_r) < 2^{-n}$, then 
$\frac{e}{8 D} \le d(z, J_r) <   2^{-n}$, and 
$\frac{e}{D} < 8 \cdot 2^{-n}$, so the algorithm outputs $1$. 
If $d(z, J_r)> 256 \cdot 2^{-n}$, then 
$\frac{16 e}{D}> d(z, J_r) > 256 \cdot 2^{-n}$. Hence 
$\frac{e}{D}> 16 \cdot 2^{-n}$, and the algorithm outputs $0$. 
\end{proof}

\begin{claim}
\label{cl14}
Step \ref{st7} is executed at most $N+1$ times, and if it outputs $1$, 
it is a valid answer. 
\end{claim}

\begin{proof}
This follows from the definition of the number $N$, the existence 
and computability of which has been established in claim
\ref{cl3}.
\end{proof}

\begin{claim}
\label{clst9}
Step \ref{st9} is executed at most a constant number of times 
between two executions of step  \ref{st7}.
\end{claim}

\begin{proof}
This is true because every point in $E_2$ is either in 
the  neighborhood of a parabolic point, or a pre-image of 
order at most $u$ of a parabolic point. Hence we can have 
at most $u$ iterations of step \ref{st9} before an iteration 
with step \ref{st7} or \ref{st11} being executed. A series of 
step \ref{st11} iterations ends with a termination or with 
a step \ref{st7} iteration before another step \ref{st9} 
iteration. 
\end{proof}

\begin{claim}
\label{clD}
After $j$ iterations, out of which $i$ are of step \ref{st7} or
\ref{st9}, $|D|>d^i$. In particular, by claims \ref{cl14} and \ref{clst9}, we always have
$|D|>d^{O(N)} = d^{O(n)}=2^{-O(n)}$. 
\end{claim}

\begin{proof}
By property \ref{pr6}, each step \ref{st7} and \ref{st9} contribute at most 
$d$ to $D$. Steps \ref{st6}, \ref{st8} and \ref{st10} terminate the 
algorithm. In step \ref{st11}, we do the long iteration only until
the angle between a repelling direction and $w$ exceeds $\al$. Until 
that moment, by property \ref{pr5n} in section \ref{nonunif}, $|r'(w)|\ge 1$, 
hence step \ref{st11} does not decrease $D$. 
\end{proof}

\begin{claim}
Step \ref{st8} outputs a valid answer for some constant $\be$, for sufficiently large $n$. 
\end{claim}

\begin{proof}
This is clearly true if $w=z$ it the $0$-th iteration. Otherwise, 
it is obvious that the previous step could not have been a step \ref{st11}, 
hence it must have been either a step \ref{st7} or \ref{st9}. In either case, 
during the previous iteration the value of $w$ was some $w'$ such that 
$r(w')=w$. By claim \ref{clinU} $w' \in r(U)$.

Denote the preparabolic point near $w$ by $q$. Then $w'$ is 
a preparabolic point $q'$ of order at most $u+1$. Since the 
parabolic points are not postcritical, there is some $\ga$ 
such that a $2 \ga$ neighborhood of any preparabolic point
of order $\le u+1$ is mapped by $r$ in a one-to-one fashion 
with no critical points. $|r'(q')| > d$, since $q' \in J_r \subset r(U)$, and 
by Koebe's theorem, $r(B(q',\ga))$ contains a ball of radius 
$d \cdot \ga/4$ around $q$. We  can take $\be$ sufficiently
large, so that $2^{-\be n} < d \cdot \ga/4$ for all $n$.

Suppose the algorithm exits on step \ref{st8}. 
Denote $d'=d(w,J_r)\le |w-q|<2^{-\be n}<d \cdot \ga/4$.  
Consider the one-to-one restriction $\ti{r}$ of $r$ to 
$B(q', 2 \ga)$. $w$ is in the $d \cdot \ga/4$-neighborhood 
of $q$, and hence $w'$ is in the $\ga$-neighborhood of $q'$. 
Denote $e=d(w', J_r)$. $e<\ga$, since $q'\in J_r$. Consider
the set $\ti{r}(B(w', e/2))$. It contains no points of 
$J_r$, and by Koebe's theorem it contains the ball $B(w, e |r'(w')|/8)$. 
Hence $d(w, J_r) \ge e |r'(w')|/8 > e \cdot d/8$. On the other hand
$d(w, J_r) < 2^{-\be n}$. Combining these inequalities we obtain 
$2^{-\be n}> e \cdot d/8$, and hence $e < 8 \cdot 2^{-\be n}/d$. 

By claim \ref{clD}, we always have $D>2^{-\eta n}$ for some constant
$\eta$. $w' \in r(U)$, and by lemma \ref{Kcol} with $w'$ we have 
$$
d(z,J_r) \le \frac{16 e}{D} < \frac{128 \cdot 2^{-\be n}/d}{2^{-\eta n}} = 
2^{7-\be n - \log d + \eta n} < 2^{-n},
$$
if we take $\be > 8 - \log d + \eta$. This is the value of $\be$ we should 
take. 
\end{proof}

\begin{claim}
Step \ref{st11} is executed at most $O(n)$
number of times between two executions of step \ref{st7}. 
\end{claim}

\begin{proof}
According to lemma \ref{iterz:main} (the long iteration lemma), we can 
make a step $\Omega(m^r)$ iterations if $|z|<\frac{1}{m}$. For a sufficiently
small $\al$, and close enough to the parabolic point,  if $|z|>\frac{1}{2 m}$,
 we have 
 $$
 |r^{m^r/C} (z)|= |z|+ \Omega(m^r |z|^{r+1}/C)=|z| + \Omega\left(\frac{m^r}{2^{r+1}
 C m^{r+1}}\right) = |z|+ \Omega\left(\frac{1}{2^{r+1} C m}\right) > |z|\cdot (1+\de), 
 $$
 for some $\de>0$. If the algorithm did not terminate at step \ref{st8}, 
 $|z|> 2^{-\be n}$, and it will take $O(n)$ long iterations to escape the
  neighborhood of the parabolic point and either terminate or reach a step \ref{st7}. 
\end{proof}

\medskip

It follows from the claims that
 
\begin{enumerate}
	\item 
	The algorithm terminates after $O(n^2)$ iterations.
	\item 
	When it terminates, it outputs a valid answer. 
\end{enumerate}

This shows that the algorithm is polynomial and correct.

\subsection{Uniformizing the construction}
\label{unif}

In this section we show how to uniformize the construction. In 
other words, we are trying to construct one machine computing
$J_{r}$ for the biggest possible family of parabolic $r$'s. 
As has been mentioned in section \ref{compjuls},
the output of the machine varies continuously in the Hausdorff metric 
with the input coefficients. The map $J: c \mapsto J_{z^2+c}$ is 
discontinuous at the parabolic point $c=1/4$ (see \cite{Douady}). 
Thus, we cannot expect one machine to compute {\em all} hyperbolic 
and parabolic sets even in the quadratic case. 

Despite the big number of different parameters that were mentioned 
as pre-requisites in section \ref{nonunif}, we will argue that 
all the information can be derived from some basic information 
about the number and periods of the parabolic points. 

First, we prove the following. 

\begin{claim}
Given the information on the parabolic orbits, we can extract the 
information on the attracting orbits ourselves. 
\end{claim}

\begin{proof}
The immediate basin of each attracting periodic orbit contains 
at least one critical point (e.g. Theorem $8.6$ in \cite{Milnor}). 
On the other hand, in our case every critical point converges either to an 
attracting or to a parabolic orbit. We proceed as follows. 
Iterate each critical point until we know to which orbit it 
converges. If it converges to an attracting orbit, we will 
eventually know it. Continue this process until the convergence
of all critical points is accounted for. 
\end{proof}

Probably the most interesting part is computing the set $U$. So 
far we haven't even shown that such a $U$ exists. To compute $U$, 
we start with a set $\ti{U}$ defined as follows. Around each 
attracting orbit we take a small ball in the basin of attraction. 
Denote the union of these balls by $\ti{A}$. 

\begin{figure}[!h]
\begin{center}
\includegraphics[angle=0, scale=0.6]{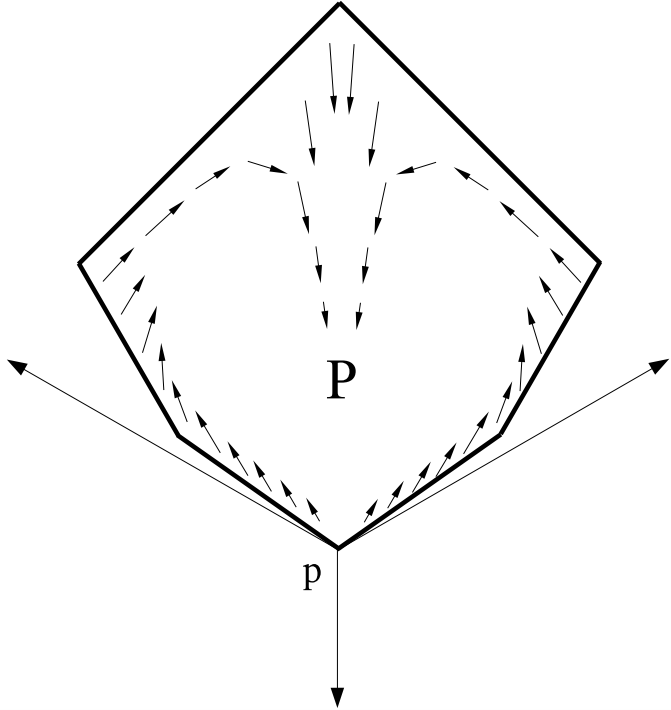}
\end{center}
\caption{The region $P$, here there are $3$ repelling directions}
\label{diamond}
\end{figure}

Around each parabolic point $p$, for any attracting 
direction $d$, consider a  small ``diamond" shaped region $P$
around $d$ such that: 
\begin{itemize}
\item 
$r(P) \subset P$ with $P \cap r(P) = \{p\}$, and 
\item
the angle of $P$ at $p$ is at least $\frac{63}{64}$ of the 
angle between the two repelling directions.
\end{itemize}
The edges of $P$ are chosen so that points would not escape 
it under $r$. 
See figure \ref{diamond} for an illustration. This is possible by the 
basic properties of the series expansion of $r$ near $p$. 

Denote the union of these ``diamons" $P$  by $\ti{P}$. Define 
$$
\ti{U} = \hat{C} - ( \ti{A} \cup \ti{P}).
$$
We know that there is an iteration $q$ such that 
\begin{enumerate}
\item 
for all critical points $c$, $r^{q-1}(c)$ is outside $\ti{U}$, and
\item 
$r^{q-1}(\infty)$ is outside $\ti{U}$.
\end{enumerate}
This is true, since the orbits of 
all the critical points eventually converge either
to an attracting or a parabolic orbit, and by our assumption 
so does the orbit of $\infty$. We have the following claims.

\begin{claim}
$\ti{V} = r^{-1} (\ti{U}) \subset \ti{U}$, and 
$\partial \ti{U} \cap \partial \ti{V} = \{ \mbox{the set of parabolic points} \}$. 
\end{claim}

\begin{proof}
This follows immediately from the definition of $\ti{U}$. 
\end{proof}

\begin{claim}
\beq
\label{unif:e1}
\bigcap_{n=0}^{\infty} r^{-n} (\ti{U}) = J_r.
\eeq
\end{claim}

\begin{proof}
The orbit of any $z\in J_r$ always stays in $J_r \subset \ti{U}$, 
hence such a $z$ is in the intersection above. 

The orbit of any $w \notin J_r$  eventually converges to 
either an attracting or a parabolic orbit, and thus escapes 
$\ti{U}$. This means that $r^k (w) \notin \ti{U}$ for some $k$, and 
$w \notin r^{-k} (\ti{U})$. So $w$ is not in the intersection in 
this case. 
\end{proof}

\begin{claim}
\label{satpr3}
For sufficiently large $q$, $U = r^{-q}(\ti{U})$ satisfies the 
conditions of part \ref{pr3} in section \ref{nonunif}.
\end{claim}

\begin{proof}
The first five conditions are satisfied automatically by the 
definition of $\ti{U}$. The last condition follows from the 
fact that for any $q$, $\partial (r^{-q}(\ti{U})) \cap \partial (r^{-q}(\ti{V}))$
consists of the parabolic points and their pre-images up to order 
$q$. 

The hardest condition to satisfy is the sixth one. Namely, we want 
to have for any $u \in r(U)$, 
	$d(u, P_r) \ge 32 \cdot d(u, J_r)$. Here $P_r$ denotes the 
	postcritical set of $r$. 
	
	Let $N_\ve$ be an $\ve$-neighborhood of 
	the parabolic points, and $M_\ve$ -- an $\ve$-neighborhood of the 
	attracting orbit points. We know that for any $\ve$ only finitely 
	many points of $P_r$ lie outside $N_\ve \cup M_\ve$. 
	For a sufficiently 
	small $\ve$ all the points in $P_r \cap N_\ve$ lie in a small 
	angular neighborhood of the attracting directions. 
	
	Denote by $\ga$ the angle between two adjacent repelling directions. 
	By the definition of $\ti{P}$,  for any $q$, all the points in $r^{-q} (\ti{U})$ 
	are in a $\frac{1}{128}\ga$-neigborhood of the attracting direction. 
	Thus the condition is satisfied in $N_\ve$. 
	
	Outside of $N_\ve \cup M_\ve$, there are only finitely many points of $P_r$, hence there is 
	a minimum $d$ of their distances from $J_r$. This minimum also exists for points in $M_\ve$, 
	since attracting orbits are bounded away from $J_r$. 
	 By \eref{unif:e1} and compactness, 
	for a sufficiently large $q$, $r^{-q}(\ti{U})$ is in the $\frac{d}{32}$-neighborhood 
	of $J_r$, and the condition is satisfied outside of $N_\ve$. 
\end{proof}

\begin{claim}
\label{cBBY}
For any constant $c$, we can produce a $2^{-c}$-precise image of $J_r$. 
\end{claim}

\begin{proof}
 This can be done using a procedure described in \cite{BBY}, Theorem $1.2$.
 Note that since $c$ does not depend on $n$, the running time of this 
 procedure will not depend on $n$ as well. 
 \end{proof}
 
 In particular, claim \ref{cBBY} immediately allows us to 
 establish property \ref{pr8} in section \ref{nonunif}.

\begin{claim}
The $q$ (and hence $U$) from claim \ref{satpr3} can be computed from the basic information 
about the parabolic points. 
\end{claim}

\begin{proof}
The proof of claim \ref{satpr3} is constructive, except for the argument 
outside of $N_\ve$, which uses compactness. There are only finitely many 
points of $P_r$ outside of $N_\ve \cup M_\ve$, which can be easily computed.
The distance from $J_r$ to the points of $M_\ve$ can also be easily 
bounded from below.  
We can find the desired value of $q$ by computing a sufficiently good
approximation of $J_r$, which is done by claim \ref{cBBY}. The precision 
with which we will have to perform this computation depends on $r$ but not  
on $n$. 
\end{proof}
\medskip

The other parts of the construction are easily seen to be uniformizable.

 \end{document}